\documentclass[10pt]{amsart}
\usepackage{amsfonts}
\usepackage[]{amsmath}
\usepackage[]{amssymb}
\usepackage[]{amsthm}	
\usepackage[]{amsfonts}

\usepackage{setspace}

\newtheorem {theorem}    {Theorem} [section]
\newtheorem {lemma}      {Lemma} [section]
\newtheorem {definition} {Definition} [section]
\newtheorem {proposition}{Proposition} [section]
\newtheorem {corollary}  {Corollary} [section]
\theoremstyle {definition} \newtheorem {remark}{Remark} [section]

\newtheorem {conjecture} {Conjecture} 
\newtheorem* {theorem*}   {Theorem} 
\newtheorem* {lemma*}      {Lemma} 
\newtheorem* {acknowledgements*} {Acknowledgements} 

\newcommand \Tan {\operatorname{Tan}}
\newcommand \Var {\operatorname{Var}}
\newcommand {\Ric} {\operatorname{Ric}}
\newcommand {\graph} {\operatorname{graph}}
\newcommand \Rm{\operatorname{Rm}}

\newcommand \inj{\operatorname{inj}}
\newcommand \sing{\operatorname{sing}}
\newcommand \reg{\operatorname{reg}}
\newcommand \spt{\operatorname{spt}}

\newcommand \tr{\operatorname{tr}}

\newcommand \dist{\operatorname{dist}}
\newcommand \mc{\operatorname{H}}
\newcommand \D{\operatorname{D}}
\newcommand \h{\operatorname{h}}
\newcommand \interior {\operatorname{int}}

\title [Generalized Apparent Horizons] {Existence, Regularity, and Properties of Generalized Apparent Horizons}
\author{Michael Eichmair}

\address{Michael Eichmair, Department of Mathematics, Massachusetts Institute of Technology, Cambridge MA-02139, USA}
\email{eichmair@math.mit.edu}

\begin{document}

\begin{abstract} We prove a conjecture of Tom Ilmanen's and Hubert Bray's on the existence of the outermost \emph{generalized apparent horizon} in an  initial data set and that it is outer area minimizing. 
\end{abstract}

\maketitle

We use the features of the construction in \cite{Plateau} to prove the following conjecture in \cite{BrayKhuri}: 

\begin{conjecture} [Hubert Bray, Tom Ilmanen, 2006] \label{conj: Ilmanen} Let $(M^n, g, p)$ be a complete asymptotically flat initial data set which contains a generalized trapped surface, and let $2 \leq n \leq 7$. Then there exists a unique outermost generalized trapped surface $\Sigma^{n-1} \subset M^n$. Moreover, this $\Sigma^{n-1}$ is a generalized apparent horizon and it is outer area minimizing (every hypersurface which encloses it has at least its area). 
\end{conjecture}


\section {Introduction, Overview, Notation}

In this paper, we consider initial data sets of general relativity, by which we will mean triples $(M^n, g, p)$ consisting of a complete oriented  Riemannian manifold $(M^n, g)$, whose regularity is at least $\mathcal{C}^3$, together with a symmetric $(0, 2)$-tensor $p$ that is required to be $\mathcal{C}^{1, \gamma}$ for some $\gamma \in (0, 1)$. We will always assume that $n \geq 2$. An initial data set is said to be asymptotically flat in the literature if the complement of a compact set in $M^n$ consists of a disjoint union of Euclidean ends $\{N_1, \ldots, N_p \}$, each diffeomorphic to $\mathbb{R}^n \setminus B(0, 1)$ with appropriate decay of the metric tensor $g_{ij}$ to the Euclidean metric $\delta_{ij}$ and the second fundamental form $p_{ij}$ to $0$ in these coordinate charts. For the purpose of Conjecture \ref{conj: Ilmanen}, the following weak decay conditions are sufficient: 
\begin{eqnarray} \label{eqn: decay asymptotically flat}
|g_{ij} - \delta_{ij}| + |x||\partial_k g_{ij}| &=& O (|x|^{-p}) \text{ and} \\ \nonumber
|p_{ij}| &=& O(|x|^{-q})
\end{eqnarray}
for some powers $p > 0$ and $q > 1$. \\

\begin{definition} [\cite{BrayKhuri}]
The compact embedded $\mathcal{C}^2$-hypersurface $\Sigma^{n-1}$ in an asymptotically flat initial data set $(M^n,g,p)$ is a generalized trapped surface with respect to a chosen end $N_k$, if it is the boundary of an open set $U \subset M^n$ which contains the `points at infinity' of all asymptotically flat ends but the chosen one, and if 
\begin{equation*}
   \mc_{\Sigma} \leq |\tr_{\Sigma}(p)|.
\end{equation*}
If $\Sigma^{n-1}$ satisfies
\begin{equation*}
   \mc_{\Sigma} = |\tr_{\Sigma}(p)|
\end{equation*}
then $\Sigma^{n-1}$ is called a generalized apparent horizon. Here, the mean curvature is computed with respect to the unit normal `pointing to infinity' so that when $U = \{x : |x| > R\} \subset \mathbb{R}^n$ one has  $\mc_\Sigma = (n-1)/R$. 
\end{definition}

\noindent The decay conditions in (\ref{eqn: decay asymptotically flat}) are chosen so that if $U$ is the complement of a large coordinate ball in an asymptotically flat manifold, then its boundary is \emph{generalized untrapped} in the sense that $\mc_{\partial U} > |\tr_{\partial U} (p)|$. Note also that by virtue of their definition, there is a partial ordering on the collection of all generalized trapped surfaces with respect to a chosen end $N_k$: if $\Sigma' = \partial U'$ and $\Sigma'' = \partial U''$ then $\Sigma''$ lies outside $\Sigma'$ if $U''$ includes $U'$. The term ``outermost" in Conjecture \ref{conj: Ilmanen} is understood in this sense. \\

Conjecture \ref{conj: Ilmanen} is well-known to be true when $p \equiv 0$, see \cite[\S4]{HuiIlm} and the references therein. In this case, generalized apparent horizons are minimal surfaces which puts variational methods at one's disposal. However, unless $p$ is a constant multiple of the metric, (generalized) apparent horizons are not known to arise as critical points of an elliptic variational problem, which complicates their existence and regularity theory. In their celebrated proof of the spacetime positive mass theorem \cite{PMTII}, R. Schoen and S.-T. Yau observed that apparent horizons can appear as an obstruction to proving the existence of entire solutions of Jang's equation on asymptotically flat initial data sets. Schoen proposed that this `defect' of Jang's equation can be turned into a device to prove existence of apparent horizons. Using this observation, L. Andersson and J. Metzger \cite{AM07} proved that closed apparent horizons exist between inner and outer trapped surfaces. The Plateau problem for apparent horizons was settled in \cite{Plateau}. The method in \cite{Plateau} shows that the apparent horizons appearing in the construction have a certain almost minimizing property (in the sense of Almgren) which likens them with minimal and constant mean curvature surfaces. Explicit bounds for their hypersurface measure in all dimensions and strong geometric estimates are a direct consequence of this property. Using very different techniques, a lower bound on the ``outer injectivity radius" for a certain class of closed apparent horizons was derived in \cite{AM07} by a delicate surgery procedure based on a priori curvature estimates. Such curvature estimates were obtained by the authors from stability and the Gauss-Bonnet theorem in \cite{AM05}, and then used for the surgery procedure to derive area bounds for certain $2$-dimensional horizons. These area bounds and the estimate on the outer injectivity radius were then applied in \cite{AM07} to show that the boundary of the trapped region of a $3$-dimensional initial data set is smooth and embedded. This important result of L. Andersson and J. Metzger is the analogue of Theorem \ref{thm: Ilmanen} for marginally outer trapped surfaces in dimension $n=3$. In this paper we adapt the features of the construction in \cite{Plateau} to generalized apparent horizons to prove Conjecture \ref{conj: Ilmanen}. Our methods here work in all dimensions $n\leq7$ and can be used to extend the result on the trapped region in \cite{AM07} to these dimensions, and without recourse to their surgery procedure. A variant of the standard calibration argument used in \cite[Example A.1]{Plateau} to establish the almost minimizing property shows that the outermost generalized apparent horizons is indeed outer area minimizing, as was conjectured by Bray and Ilmanen. \\

Our basic existence result for generalized apparent horizons is as follows:

\begin{theorem} \label{thm: basic existence introduction} Assume that $2 \leq n \leq 7$ and let $(M^n, g, p)$ be a complete initial data set. Let $\Omega \subset M^n$ be a bounded open subset with embedded $\mathcal{C}^2$-boundary $\partial \Omega = \partial_1 \Omega \dot{\cup} \partial_2 \Omega$ such that $\mc_{\partial_1 \Omega} > |\tr_{\partial_1 \Omega}(p)|$ (mean curvature with respect to the normal pointing out of $\Omega$) and $\mc_{\partial_2 \Omega} \leq |\tr_{\partial_2 \Omega}(p)|$ (with respect to the normal pointing  into $\Omega$). Then there exists an open set $\Omega' \subset \Omega$ such that $\partial \Omega' = \partial_1 \Omega \dot \cup \Sigma^{n-1}$ where $\Sigma^{n-1}$ is a closed embedded $\mathcal{C}^{2, \alpha}$-hypersurface satisfying $\mc_{\Sigma} = |\tr_{\Sigma}(p)|$ (normal pointing into $\Omega'$) with uniform estimates depending only on $\mathcal{H}^{n-1}(\partial_1 \Omega)$, $|p|_{\mathcal{C}^1(\bar\Omega)}$, and the local geometry of $(M^n, g)$. If a connected component of $\Sigma^{n-1}$ intersects with a component of $\partial_2 \Omega$, then these components coincide. Moreover, $\Sigma^{n-1}$ minimizes area with respect to variations in $\bar {\Omega}'$ and it is $C$-almost minimizing in all of $\Omega$ for a constant $C = C(|p|_{\mathcal{C}(\bar \Omega)})$. 
\end{theorem}

\noindent The existence statement in Theorem \ref{thm: basic existence introduction} is the analogue for generalized apparent horizons of \cite[Theorem 3.1]{AM07}. Our proof of Theorem \ref{thm: basic existence introduction} in Section \ref{sec: existence of generalized apparent horizons} is based on the Perron method used in \cite{Plateau} together with an approximation argument. The $C$-almost miniziming property with the constant $C$ depending only on $|p|_{\mathcal{C}(\bar \Omega)}$ and the outer area minimizing property in the preceding theorem are important features of this approach and at the heart of our proof of Conjecture \ref{conj: Ilmanen}. A summary of results related to what we call the $C$-almost minimizing property here and the class $\mathcal{F}_C$ of boundaries with this property is given in \cite[Appendix A]{Plateau} with concise references to the geometric measure theory literature. We derive the $\mathcal{C}^{2, \alpha}$-estimates from results in \cite{SchoenSimon}, which we appropriate to our context in Appendix \ref{sec: remark on SS81}. The robust `low order approach' to regularity used here via geometric measure theory and the stability based analysis of  \cite{SchoenSimon} is available also when $n \geq 8$ if we accept thin singular sets. It provides a satisfactory theory for limits of regular embedded horizons in arbitrary dimensions that is friendly towards analysis, see Remark \ref{rem: curvature estimates follow too} for details, and compare with the curvature estimates that were obtained in \cite{AM05} by generalizing the iteration method of \cite{Schoen-Simon-Yau}. The modification of Jang's equation used in this paper can be applied to find solutions to the Plateau problem for generalized apparent horizons as in \cite{Plateau}. \\

\noindent The outer area minimizing property of generalized apparent horizons is appealing from the point of view of a spacetime Penrose inequality, and we refer the reader to the paper \cite{BrayKhuri} for their original motivation in this context. See also the recent construction by A. Carrasco and M. Mars in \cite{Mars} of a counterexample to a conjectured spacetime Penrose inequality in \cite{BrayKhuri}. It has been shown by M. Khuri in \cite{K09} that $3+1$-dimensional Minkowski space does not contain generalized apparent horizons. G. Galloway has an argument to construct generalized trapped \emph{curves} in $2+1$-dimensional Minkowski space (private communication). \\

In Section \ref{sec: existence of outermost generalized apparent horizon} we observe that the limit of an increasing sequence of generalized apparent horizons $\{\Sigma^{n-1}_i\}$, where the $\Sigma^{n-1}_i$ satisfy uniform $\mathcal{C}^{2, \alpha}$-estimates and are all outer area minimizing, is embedded. Embeddedness is a critical issue here because (generalized) apparent horizons, unlike minimal surfaces, do not satisfy a two-sided maximum principle, cf. \cite[Remark 8.3]{AM05}. This observation leads directly to the proof of Conjecture \ref{conj: Ilmanen}: 

\begin{theorem} \label{thm: Ilmanen} Conditions as in the preceding theorem. Then there exists a unique outermost generalized apparent horizon $\Sigma^{n-1}$ in $\Omega$. $\Sigma^{n-1}$ is a closed embedded $\mathcal{C}^{2, \alpha}$-hypersurface. Moreover, $\Sigma^{n-1}$ is outer area minimizing.  
\end{theorem}

\noindent A refined statement and the proof of this theorem are given in Section \ref{sec: existence of outermost generalized apparent horizon}. \\

The following notation will be used throughout this paper:

\begin{definition} Given an initial data set $(M^n, g, p)$, an open set $U \subset M^n$, and a $\mathcal{C}^2$-function $u : U \to \mathbb{R}$ we write 
\begin{eqnarray*}
\mc(u) := \D_i \left(\frac{\D^i u}{ \sqrt{1 + |\D u|^2}} \right) \\
\tr(p)(u) := \left(g^{i j}  - \frac{\D^i u \D^j u}{1 + |\D u|^2}\right) p_{ij}
\end{eqnarray*} and for every $\varepsilon > 0$
\begin{equation*}
|\tr(p)|_\varepsilon (u) := \sqrt{\tr(p)(u)^2 + \varepsilon^2}. 
\end{equation*}
\end{definition}

\noindent Note that $\mc(u)$ is the mean curvature of $\graph(u, U) \subset M^n \times \mathbb{R}$ computed with respect to the product metric $g + dx_{n+1}^2$ and its downward pointing unit normal $\frac{(\D u, -1)}{\sqrt{1 + |\D u|^2}}$. If we think of the tensor $p$ as being extended to the product manifold $M^n \times \mathbb{R}$ by zero in the vertical direction, then $\tr(p)(u)$ represents its trace over the tangent space of $\graph(u, U)$. We introduce the auxiliary expressions $|\tr(p)(u)|_\varepsilon$, which regularize the absolute value of $\tr(p)(u)$, to facilitate exposition and analysis in the paper. \\

\begin{acknowledgements*} This work forms part of my thesis, and I am very much indebted to my adviser Richard Schoen for his constant support and encouragement. Thanks very much to Hubert Bray for great discussions, and for drawing my attention to this problem. I would like to sincerely thank Simon Brendle, Leon Simon, and Brian White for the great example they set and for everything they have taught me over the years. I am grateful to Greg Galloway, Tom Ilmanen, Marcus Khuri, and Jan Metzger for their interest in this work. 
\end{acknowledgements*}


\section {Solving $(\mc - |\tr(p)|_\varepsilon - t) (u) = 0$ for $t > 0$} \label{sec: existence of u_t^epsilon}

In this section we use the Perron method as in \cite{Plateau} to prove the following existence result:

\begin{lemma} \label{lem: existence of u_t^epsilon} Let $(M^n, g, p)$ be a complete initial data set and let $\Omega \subset M^n$ be a bounded open subset with embedded $\mathcal{C}^2$-boundary such that $\partial \Omega = \partial_1 \Omega \dot{\cup} \partial_2 \Omega$ where $\mc_{\partial_1 \Omega} > |\tr_{\partial_1 \Omega}(p)|_\varepsilon$ (with respect to the outward unit normal) and $\mc_{\partial_2 \Omega} < |\tr_{\partial_2 \Omega}(p)|_\varepsilon$ (with respect to the normal into $\Omega$) for some $\varepsilon \in (0,1)$. There exists $\theta > 0$ such that for every sufficiently small $t>0$ there is a non-positive function $u_t^\varepsilon \in \mathcal{C}^{2, \alpha}_{loc}(\Omega) \cap \mathcal{C}^{0, 1}(\Omega \cup \partial_1 \Omega)$ satisfying $\left(\mc - |\tr(p)|_\varepsilon- t\right) u_t^\varepsilon = 0$ in $\Omega$ with $u_t^\varepsilon = 0$ along $\partial_1\Omega$, $u_t^\varepsilon \leq - \frac {\theta}{t}$ on $\{x \in \Omega : \dist(x, \partial_2 \Omega) \leq \theta\}$, $u_t^\varepsilon \geq \ln\left(1 - \frac{\dist(\partial_1 \Omega, \cdot)}{\theta}\right)$ on $\{ x \in \bar{\Omega} : \dist(x, \partial_1 \Omega) \leq \theta\}$, and $0 \geq u^\varepsilon_t \geq - \frac{C}{t}$ on all of $\Omega$ where $C := 1 + n|p|_{\mathcal{C}(\bar \Omega)}$.  
\end{lemma}
\noindent We discuss the required modifications of the method used in \cite[\S 2, \S 3]{Plateau} for the proof of Lemma \ref{lem: existence of u_t^epsilon}. We refer the reader to that paper for details and concise references for the techniques involved. Recall first that if $u : U \subset \Omega \to \mathbb{R}$ is a $\mathcal{C}^3$-function and if $G^n := \graph(u, U) = \{(x, u(x)) : x \in U \} \subset M^n \times \mathbb{R}$ (with the product metric $g + dx_{n+1}^2$) denotes its graph, then
\begin{equation*}
\Delta_{G} \frac{1}{v} + \left (|\h|^2 + \Ric_{M \times \mathbb{R}} (\vec{\nu}, \vec{\nu}) + \vec{\nu} \mc(u) \right)  \frac{1}{v}= 0
\end{equation*}
where $v := \sqrt{1 + |\D u|^2}$ is the `area stretch factor' of the graph, cf. \cite[Equation (2.28)]{PMTII}. This is the second variation formula of the area element of $G^n$ applied to the variation coming from vertical translation in  $M^n \times \mathbb{R}$. Here, $\Delta_{G}$ denotes the non-positive hypersurface Laplacian with respect to the induced metric $\bar{g}$ on $G^n$, $\h$ denotes the second fundamental form of $G^n$, and the downward pointing unit normal $\vec{\nu} = (1 + |\D u|^2)^{-\frac{1}{2}} (\D u, -1)$ differentiates the mean curvature $\mc(u)$ of $G^n$ as a function on $U \times \mathbb{R}$ (invariant in the vertical direction). If $t > 0$ and if $u$ solves $\left(\mc - |\tr(p)|_\varepsilon-t\right)u = 0$ on $U$, then this identity implies a differential inequality of the form
\begin{eqnarray} \label{inequ: area strech factor is subharmonic}
\left( 1 - \frac{1}{3(n-1)}\right) \frac{|\h|^2}{ v}+\Delta_{G} \frac{1}{v} \leq \frac{\kappa^2}{v}
\end{eqnarray} for $v^{-1}$, where $\kappa$ is a constant depending only on $|\Ric_M|$, $|p|_{\mathcal{C}^1}$, and the dimension $n$, but \emph{not} on $\varepsilon > 0$ or $t > 0$. The scale of $\kappa$ is one over length. The method of Korevaar-Simon (see \cite[Lemma 2.1]{Plateau} for a precise statement) provides interior gradient estimates from oscillation bounds for such functions $u$. \\

The pointwise inequality (\ref{inequ: area strech factor is subharmonic}) implies a stability-type inequality for $G^n$ in a standard way: given a test function $\phi \in \mathcal{C}^1_c (U \times \mathbb{R})$, multiply (\ref{inequ: area strech factor is subharmonic}) by $v \phi^2 $ and integrate by parts on $G^n$ to obtain
\begin{equation} \label{inequ: stability of graph(u)}
\left(1 - \frac{1}{3(n-1)}\right) \int_G |\h|^2 \phi^2 \leq \int_G |\bar{\D} \phi|^2  + \kappa^2 \int_G \phi^2
\end{equation}
adjusting $\kappa$ depending on the dimension. A related computation is implicit in the proof of \cite[Proposition 1]{PMTII}, see also \cite[Lemma 5.6]{AM05}. Instead of $\frac{1}{3(n-1)}$ we could take any other positive constant less than $\frac{1}{(n-1)}$ to make the arguments below work; this cut-off however is critical for the use of \cite{SchoenSimon}, as we discuss in the appendix. \\

We now construct the particular solutions $u_t^\varepsilon$ described in Lemma \ref{lem: existence of u_t^epsilon}. First note that the constant function $0$ is a super solution of $\left( \mc - |\tr(p)|_\varepsilon - t \right) u = 0$. Similarly, the constant function at height $-\frac{C}{t}$ where $C= 1 + n|p|_{\mathcal{C}(\bar \Omega)}$ (the dimension times the largest eigenvalue of $p$ on $\bar \Omega$ plus one) is a sub solution. Let $\dist(\partial_2 \Omega, \cdot)$: $ M^n \to \mathbb{R}$ measure geodesic distance in $(M^n, g)$ to $\partial_2 \Omega$ and choose $\delta > 0$ so small that $\dist(\partial_2 \Omega, \cdot)$ is a $\mathcal{C}^2$-function on $\{x \in \Omega : \dist(\partial_2 \Omega, x) \leq 2 \delta \}$ (cf. \cite[Lemma 14.16]{GT}), such that the latter set is disjoint from $\partial_1 \Omega$, and such that the mean curvature of the distance surfaces $(\partial_2 \Omega)_{\gamma} := \{x \in \Omega : \dist (\partial_2 \Omega, x) = \gamma \}$ satisfies $2 \delta + \mc_{(\partial_2 \Omega)_\gamma} < |\tr_{(\partial_2 \Omega)_\gamma} (p)|_\varepsilon$ for all $\gamma \in [0, 2 \delta)$. It follows that the function $\bar{u}_t \in \mathcal{C}(\bar \Omega)$ defined by 
\begin{equation*}
\bar u_t^\varepsilon (x) := \left\{ \begin{array} {rl} \frac{\dist(\partial_2 \Omega, x) - \delta}{t} & \text  { if } d_{\partial_2 \Omega}(x) \leq \delta  \\ 0 & \text{ if } d_{\partial_2 \Omega}(x) > \delta \end{array}  \right.
\end{equation*}
is a Perron super solution for the operator $L_t^\varepsilon = \mc - |\tr(p)|_\varepsilon - t$ on $\Omega$ for $t>0$ sufficiently small. Similarly, it follows that for all $t>0$ sufficiently small the function defined by 
\begin{equation*}
\underline u_t (x) := \left\{ \begin{array} {rl} \ln\left(1 - \frac{\dist(\partial_1 \Omega, x)}{\delta}\right) & \text  { if } 0 \leq d_{\partial_1 \Omega}(x) \leq \delta\left(1 - \exp(-\frac{C}{t})\right)  \\ -\frac{C}{t} & \text{ if } d_{\partial_1 \Omega}(x) > \delta\left(1 - \exp(-\frac{C}{t})\right) \end{array}  \right.
\end{equation*}
is a Perron sub solution (possibly for some smaller $\delta>0$ depending on the geometry of $\partial_1 \Omega$). We emphasize that this sub solution is independent of $\varepsilon > 0$ and also of $t > 0$, when sufficiently small, in a fixed neighborhood of $\partial_1 \Omega$. Also observe that $\bar u_t^\varepsilon \leq - \frac{\delta}{2t}$ on $\{x \in \bar \Omega : 0 \leq \dist(\partial_2 \Omega, x) \leq \frac{\delta}{2} \}$, and that $u_t^\varepsilon \geq \ln\left(1 - \frac{\dist(\partial_1 \Omega, \cdot)}{\delta}\right)$ on $\{ x \in \bar{\Omega} : 0 \leq \dist(\partial_1 \Omega, x) \leq \frac{\delta}{2}\}$. The point in the construction of both sub and super solution is that the non-constant portions of their graphs converge to the cylinders $(\partial_i \Omega)_\gamma \times \mathbb{R}$ uniformly in $\mathcal{C}^2$ on compact sets, as the parameters $\delta, t$ tend to zero. This classical construction of boundary barriers for the prescribed mean curvature equations is due to J. Serrin \cite{Serrin69}, see \cite[\S14]{GT}, cf. \cite[Proposition 3.5]{AM07}, and also \cite[Lemmas 2.2 and 3.3]{Plateau} for concise references and a `geometric' discussion of these barriers in this context. \\

\noindent The Perron method employed in \cite{Plateau} now carries over verbatim to the present context and shows that there exists a $\mathcal{C}^{2, \alpha}_{loc}$-function $u_t^\varepsilon : \Omega \to \mathbb{R}$ with $\underline u_t \leq u_t^\varepsilon \leq \bar u_t^\varepsilon$ and such that $\left(\mc - |\tr(p)|_\varepsilon - t\right) u_t^\varepsilon = 0$ on $\Omega$. It follows easily that $u_t$ extends to a $\mathcal{C}^{0, 1}$ function near $\partial_1 \Omega$. Taking $\theta = \delta/2$ this concludes the proof of Lemma \ref{lem: existence of u_t^epsilon}. \\

As in \cite{PMTII} we will eventually pass the graphs of solutions of the equation $\left(\mc - |\tr(p)|_\varepsilon - t\right)u_t^\varepsilon = 0$ to a geometric subsequential limit as $t \searrow 0$. In order to avoid a dimensional restriction, we take this limit in a certain class of almost minimizing boundaries (in the sense of Almgren) and use compactness and regularity results from geometric measure theory to analyze the limit. More precisely, we use the classes $\mathcal{F}_C$ of $C$-almost minimizing boundaries discussed in \cite[Appendix A]{Plateau}. As we will see below, we are able to bound $|tu^\varepsilon_t|$ and hence the mean curvature of all graphs in the construction  uniformly by a constant that only depends on $|p|_{\mathcal{C}(\bar \Omega)}$ so that their graphs are $C$-almost minimizing, cf. \cite[Example A.1]{Plateau}. Exactly as in the analysis of \cite[Proposition 4]{PMTII}, the limiting surface will contain a graphical component that is asymptotic to a vertical cylinder, and the cross section $\Sigma^{n-1}_\varepsilon \subset M^n$ of this cylinder satisfies the geometric equation $\mc_{\Sigma_\varepsilon} = |\tr_{\Sigma_\varepsilon}(p)|$ and inherits the stability-type inequality (\ref{inequ: stability of graph(u)}) as well as the $C$-almost minimizing property in $\Omega$ from the original graphs. \\


\section {Existence of surfaces $\Sigma_{\varepsilon}^{n-1}$ with $\mc_{\Sigma_{\varepsilon}} = |\tr_{\Sigma_\varepsilon} (p)|_\varepsilon$} \label{sec: existence of Sigma_epsilon}

Following \cite{PMTII}, we now pass the graphs of the solutions $u_t^\varepsilon$ from Lemma \ref{lem: existence of u_t^epsilon} to a geometric limit as $t \searrow 0$. Just as with the regular Jang's equation, the vertically unbounded limiting surface will be asymptotic to  a cylinder whose cross section is a closed surface $\Sigma_\varepsilon^{n-1}$ such that $\mc_{\Sigma_{\varepsilon}} = |\tr_{\Sigma_{\varepsilon}}(p)|_\varepsilon$, see also \cite[Theorem 1.1]{AM07}, \cite[Theorem 1.1]{Plateau}. 

\begin{theorem} \label{thm: existence of Sigma_epsilon} Assume that $2 \leq n \leq 7$, let $(M^n, g, p)$ be a complete initial data set, and let $\Omega \subset M^n$ be a bounded open subset with embedded $\mathcal{C}^2$-boundary $\partial \Omega = \partial_1 \Omega \dot{\cup} \partial_2 \Omega$ where $\mc_{\partial_1 \Omega} > |\tr_{\partial_1 \Omega}(p)|_\varepsilon$ (with respect to the outward unit normal) and $\mc_{\partial_2 \Omega} < |\tr_{\partial_2 \Omega}(p)|_\varepsilon$ (with respect to the normal pointing into $\Omega$). There exists an open set $\Omega' \subset \Omega$ with $\partial \Omega' = \partial_1 \Omega \dot{\cup} \Sigma_\varepsilon^{n-1}$ such that $\Sigma_\varepsilon^{n-1}$ is an embedded hypersurface disjoint from $\partial \Omega$ that satisfies $\mc_{\Sigma_\varepsilon} = |\tr_{\Sigma_\varepsilon}(p)|_\varepsilon$ (unit normal pointing into $\Omega'$). There are $\mathcal{C}^{2, \alpha}$-estimates for the surface $\Sigma^{n-1}_\varepsilon$ arising in the construction which only depend on $\mathcal{H}^{n-1}(\partial_1 \Omega)$, $|p|_{\mathcal{C}^1(\bar \Omega)}$, and the local geometry of $(M^n, g)$. In fact, $\Sigma_\varepsilon^{n-1}$ is stable in the sense of Appendix \ref{sec: remark on SS81}, it minimizes area with respect to variations in $\bar\Omega'$, and it is $C$-almost minimizing in all of $\Omega$ for a constant $C = C(|p|_{\mathcal{C}(\bar \Omega)})$. 

\begin{proof} Let $\theta > 0$, $u_t^\varepsilon \in \mathcal{C}^{2, \alpha}_{loc}(\Omega) \cap \mathcal{C}^{0, 1}(\Omega \cup \partial_1 \Omega)$, $C>0$ be as in  Lemma \ref{lem: existence of u_t^epsilon}, and let $G_t^\varepsilon := \{(x, u_t^\varepsilon(x)) : x \in \Omega \cup \partial_1 \Omega\}$ be the corresponding graphs in $(M^n \times \mathbb{R}, g + dx_{n+1}^2)$. Since $ 0 \geq u_t^\varepsilon \geq - \frac{C}{t}$ the mean curvatures of the graphs $G_t^\varepsilon$ are bounded uniformly by $2C$ and hence are $2C$-almost minimizing in the language of \cite[Appendix A]{Plateau}. Using the compactness and regularity theory for these almost minimizing boundaries (see also \cite[Remark 4.1]{Plateau}) we can pass $G_{t}^\varepsilon$ to a smooth subsequential limit $G^\varepsilon$ along a sequence $t' \searrow 0$. Moreover, the connected components of $G^\varepsilon$ are either cylindrical or themselves entirely graphical by the Harnack principle, cf. \cite[Lemma 2.3]{Plateau} for reference. In fact we see that there exists an open subset $\Omega' \subset \Omega$ with $\partial \Omega' = \partial_1 \Omega \dot \cup \Sigma_\varepsilon^{n-1}$ so that below $\Omega' \times \{0\}$ the hypersurface $G^\varepsilon$ is given as the graph of a function $u^\varepsilon : \Omega' \cup \partial_1 \Omega \to \mathbb{R}$ satisfying $\mc(u^\varepsilon) = |\tr(p)|_\varepsilon (u^\varepsilon)$. This is because the sub solution  $\underline{u}_t^\varepsilon$ is independent of $t$ in a fixed neighborhood of $\partial_1 \Omega$ preventing the limit from diverging downwards there. We have that $u^\varepsilon = 0$ on $\partial_1 \Omega$ and that $u^\varepsilon$ tends to negative infinity with $\graph(u^\varepsilon, \Omega')$ asymptoting $\Sigma^{n-1}_\varepsilon \times \mathbb{R}$ on approach to $\Sigma^{n-1}_\varepsilon$ through $\Omega'$, cf. \cite[Proposition 4]{PMTII}. It follows also that $\Sigma^{n-1}_\varepsilon$ satisfies $\mc_{\Sigma_\varepsilon} = |\tr_{\Sigma_\varepsilon}(p)|$ where the mean curvature is computed with respect to the unit normal pointing into $\Omega'$. Using that the unit normal vector field $\left(1 + |\D u^\varepsilon|^2 \right)^{- \frac{1}{2}} (\D u^\varepsilon, -1)$ of $G^\varepsilon$ has non-negative divergence, a standard calibration argument (cf.  \cite[Example A.1]{Plateau}) shows that $\Sigma^{n-1}_\varepsilon$ minimizes area with respect to variations in $\bar \Omega'$. The remarks in Section \ref{sec: existence of u_t^epsilon} imply that $\Sigma^{n-1}_\varepsilon$ satisfies a stability-type inequality (\ref{inequ: stability of graph(u)}) with a constant $\kappa$ that only depends on $\Ric_M$ and $|p|_{\mathcal{C}^1(\bar \Omega)}$ making the results discussed in Appendix \ref{sec: remark on SS81} available. It follows that the locally defining functions of $\Sigma_\varepsilon^{n-1}$ satisfy uniform $\mathcal{C}^{1, \alpha}$-estimates; since they also satisfy the (geometric) divergence form equation $\mc_{\Sigma_\varepsilon} = |\tr_{\Sigma_\varepsilon}(p)|$, $\mathcal{C}^{2, \alpha}$-estimates follow from standard elliptic theory \cite{GT}. 
\end{proof}
\end{theorem}


\section {Existence of Generalized Apparent Horizons} \label{sec: existence of generalized apparent horizons}

In this section we prove an existence theorem for generalized apparent horizons by passing the hypersurfaces $\Sigma_\varepsilon^{n-1}$ constructed in the previous section to a subsequential limit as $\varepsilon \searrow 0$. The lower order geometric properties ($C$-almost minimizing, outer area minimizing) of the surfaces $\Sigma_\varepsilon^{n-1}$descend to this limit, as does the stability-type inequality (\ref{inequ: stability of graph(u)}) with uniform constant $\kappa = \kappa(\Ric_M, |p|_{\mathcal{C}^1(\bar \Omega)}, n)$, so that we stay in the class of surfaces to which the regularity and compactness theory of \cite{SchoenSimon} discussed in Appendix \ref{sec: remark on SS81} is applicable. In Proposition \ref{prop: to any two generalized apparent horizons that intersect there is another one enclosing them} we use the Perron method to prove that given two generalized trapped surfaces, there always exists a stable outer minimizing generalized apparent horizon enclosing both of them. The purpose of this proposition in the proof of Conjecture \ref{conj: Ilmanen} corresponds to that of Lemma 8 in \cite{KrieHay} and more specifically to that of Lemma 7.7 in \cite{AM07}.

\begin{theorem} \label{thm: existence of generalized apparent horizon in Omega} Assume that $2 \leq n \leq 7$ and let $(M^n, g, p)$ be a complete initial data set. Let $\Omega \subset M^n$ be a bounded open subset with embedded $\mathcal{C}^2$-boundary $\partial \Omega = \partial_1 \Omega \dot{\cup} \partial_2 \Omega$ such that $\mc_{\partial_1 \Omega} > |\tr_{\partial_1 \Omega}(p)|$ (mean curvature with respect to the normal pointing out of $\Omega$) and $\mc_{\partial_2 \Omega} \leq |\tr_{\partial_2 \Omega}(p)|$ (with respect to the normal pointing into $\Omega$). Then there exists an open set $\Omega' \subset \Omega$ such that $\partial \Omega' = \partial_1 \Omega \dot \cup \Sigma^{n-1}$ where $\Sigma^{n-1}$ is a closed embedded hypersurface satisfying $\mc_{\Sigma} = |\tr_{\Sigma}(p)|$ (normal pointing into $\Omega'$) with uniform $\mathcal{C}^{2, \alpha}$-estimates depending only on $\mathcal{H}^{n-1}(\partial_1 \Omega)$, $|p|_{\mathcal{C}^1(\bar \Omega)}$, and the local geometry of $(M^n, g)$. The hypersurface $\Sigma^{n-1}$ is disjoint from $\partial_1 \Omega$, and if a connected component of $\Sigma^{n-1}$ intersects with a component of $\partial_2 \Omega$, then these components must coincide. Moreover, $\Sigma^{n-1}$ minimizes area with respect to variations in $\bar {\Omega}'$ and it is $C$-almost minimizing in all of $\Omega$ for a constant $C = C(|p|_{\mathcal{C}(\bar \Omega)})$. 

\begin{proof} By the result of Section \ref{sec: existence of Sigma_epsilon} we can find for every sufficiently small $\varepsilon > 0$ an open set $\Omega'_\varepsilon \subset \Omega$ with embedded boundary $\partial \Omega'_\varepsilon = \partial_1 \Omega \dot \cup \Sigma^{n-1}_\varepsilon$ with  $\mc_{\Sigma_\varepsilon} = |\tr_{\Sigma_\varepsilon}(p)|$ (normal pointing into $\Omega'_\varepsilon$) so that (i) $\Sigma^{n-1}_\varepsilon \in  \mathcal{F}_C(\Omega)$
with a constant $C$ independent of $\varepsilon >0$, (ii) $\Sigma^{n-1}_\varepsilon$ minimizes area with respect to variations in $\bar{\Omega}'_\varepsilon$, and (iii) $\Sigma^{n-1}_\varepsilon$ is $\mathcal{C}^{2, \alpha}$ with estimates independent of $\varepsilon > 0$. \\

\noindent By property (iii) (alternatively using Theorem \ref{thm: varifold compactness for stable submanifolds}) we can pass $\Sigma_{\varepsilon}^{n-1}$ to a $\mathcal{C}^{2, \beta}$-subsequential limit $\Sigma^{n-1}$, where we fix some $\beta \in (0, \alpha)$. The $C$-almost minimizing property (i) of $\Sigma_\varepsilon^{n-1}$ gives that the only place where $\Sigma^{n-1}$ could fail to be embedded is along $\partial_2 \Omega$. However, we can use the convergence in $\mathcal{C}^{2, \beta}$ and property (i) to rule out sheeting near the boundary exactly as in the proof of Corollary \ref{cor: increasing limits are embedded}. It follows that $\Sigma^{n-1}$ is properly embedded, that $\Sigma^{n-1} = \partial \Omega'$ for an open set $\Omega' \subset \Omega$, and that $\mc_{\Sigma} = |\tr_\Sigma(p)|$. The strong maximum principle shows in a standard way that if a component of $\Sigma^{n-1}$ intersects with a component of $\partial_2 \Omega$, then these components must coincide. (See \cite[Proposition 3.1]{AM07} \cite[Proposition 2.4]{AG05} for marginally outer trapped surfaces; in the notation of the latter reference, $0 \leq u_2 - u_1$ (the difference of the two locally defining functions) satisfies a linear elliptic equation to which the strong maximum principle can be applied.) That $\Sigma^{n-1}$ minimizes area with respect to variations in $\bar{\Omega}'$ and is $C$-almost minimizing in all of $\Omega$ follow easily now. 
\end{proof}

\end{theorem}

\begin{lemma} \label{lem: kicking out the kink} Let $(M^n, g, p)$ be a complete initial data set, $2 \leq n \leq 7$, and let $\Omega', \Omega'' \subset \Omega \subset M^n$ be open bounded subsets so that $\partial \Omega = \partial_1 \Omega \dot \cup \partial_2 \Omega$ with $\mc_{\partial_1 \Omega} > |\tr_{\partial_1 \Omega}(p)|$ and such that $\partial \Omega' = \partial_1 \Omega \dot \cup \Sigma'$ and $\partial \Omega'' = \partial_1 \Omega \dot \cup \Sigma''$ where $\mc_{\Sigma'} = |\tr_{\Sigma'}(p)|$ (normal pointing into $\Omega'$) and $\mc_{\Sigma''} = |\tr_{\Sigma''}(p)|$ (normal pointing into $\Omega''$) for embedded $\mathcal{C}^{2}$-hypersurfaces $\partial_1 \Omega, \Sigma', \Sigma''$. Then for every sufficiently small $\varepsilon > 0$ there exists an open set $\Omega'''_\varepsilon \subset \Omega' \cap \Omega''$ such that $\partial \Omega'''_\varepsilon = \partial_1 \Omega \dot \cup \Sigma'''_\varepsilon$ where $\Sigma'''_\varepsilon$ is embedded and satisfies $\mc_{\Sigma'''_\varepsilon} = |\tr_{\Sigma'''_\varepsilon}(p)|_\varepsilon$ (normal pointing into $\Omega'''_\varepsilon$). The surface $\Sigma'''_\varepsilon$ arising in our construction satisfies $\mathcal{C}^{2, \alpha}$-estimates depending only on $\mathcal{H}^{n-1}(\partial_1 \Omega)$, $|p|_{\mathcal{C}^1(\bar\Omega)}$, and the local geometry of $(M^n, g)$, but not on $\varepsilon > 0$. Moreover, $\Sigma'''_\varepsilon$ minimizes area with respect to variations in $\bar{\Omega}'''_\varepsilon$, and it is $C$-almost minimizing in all of $\Omega' \cap \Omega''$ for some $C = C(|p|_{\mathcal{C}(\bar \Omega)})$.

\begin{proof} As in Section \ref{sec: existence of u_t^epsilon} we can use the Perron method to construct for every sufficiently small $t>0$ a non-positive function $u_t'^\varepsilon \in \mathcal{C}^{2, \alpha}_{loc}(\Omega') \cap \mathcal{C}^{0, 1} (\Omega' \cup \partial_1 \Omega)$ such that $(\mc - |\tr(p)|_\varepsilon - t) u_t'^\varepsilon = 0$ with $u_t'^\varepsilon = 0$ on $\partial_1 \Omega$. Moreover, we can arrange that $u_t'^\varepsilon \geq - \frac{C}{t}$ on all of $\Omega'$ and  $u_t'^\varepsilon \leq - \frac{\theta}{t}$ on $\{x \in \Omega' : \dist(x, \Sigma') < \theta \}$ where $C = C(|p|_{\mathcal{C}(\bar\Omega)})$ and where $\theta > 0$ depends on $\varepsilon >0$ but not on $t$. We construct $u_t''^\varepsilon$ with identical properties, but with respect to $\Omega''$. \\

\noindent We now consider the function $\bar{u}_t'''^\varepsilon = \min (u_t'^\varepsilon, u_t''^\varepsilon) \in \mathcal{C}(\Omega' \cap \Omega'')$. Then $\bar{u}_t'''^\varepsilon$ is a Perron super solution with respect to the operator $\mc - |\tr(p)|_\varepsilon - t$ on the open set $\Omega' \cap \Omega''$. We also have that $\bar{u}_t'''^\varepsilon \leq - \frac{\theta}{t}$ on the set $\{x \in \Omega' \cap \Omega'' : \dist(x, \Sigma') < \theta \text{ or } \dist(x, \Sigma'') < \theta \}$, so that in particular $\bar{u}_t'''^\varepsilon$ tends to negative infinity on approach to $\partial (\Omega' \cap \Omega'') \setminus \partial_1 \Omega$ uniformly as $t \searrow 0$. Let $u_t^\varepsilon$ be the Perron solution constructed from this super solution and the Perron sub solution described in Section \ref{sec: existence of u_t^epsilon} (which is constant away from $\partial_1 \Omega$ and independent of $t, \varepsilon > 0$ near $\partial_1 \Omega$). The graphs $G_t^\varepsilon := \graph(u_t^\varepsilon, \Omega' \cap \Omega'') \subset \Omega \times \mathbb{R}$ are $C$-almost minimizing in the cylinder $(\Omega' \cap \Omega'') \times \mathbb{R}$. As in the proof of Theorem \ref{thm: existence of Sigma_epsilon} we conclude that there exists an open set $\Omega'''_\varepsilon \subset \Omega' \cap \Omega''$ with $\partial \Omega'''_\varepsilon = \partial_1 \Omega \dot \cup \Sigma'''_\varepsilon$ such that $\Sigma'''_\varepsilon$ is embedded, outer area minimizing with respect to $\bar \Omega '''_\varepsilon$, $C$-almost minimizing in all of $\Omega' \cap \Omega''$, and so that $\mc_{\Sigma'''_\varepsilon} = |\tr_{\Sigma'''_\varepsilon}(p)|_\varepsilon$. The $\mathcal{C}^{2, \alpha}$-estimates follow from the comments succeeding the statement of Lemma \ref{lem: existence of u_t^epsilon} and the results in Appendix \ref{sec: remark on SS81}. \end{proof}
\end{lemma}

\begin{proposition} \label{prop: to any two generalized apparent horizons that intersect there is another one enclosing them} Assumptions as in Lemma \ref{lem: kicking out the kink}. Then there exists an open set $\Omega''' \subset \Omega' \cap \Omega''$ with embedded boundary $\partial \Omega''' = \partial_1 \Omega \dot \cup \Sigma'''$ such that $\Sigma'''$ is disjoint from the intersecting (but not coinciding) components of $\Sigma', \Sigma''$ and so that $\mc_{\Sigma'''} = |\tr_{\Sigma'''}(p)|$ (with respect to the normal pointing into $\Omega'''$). The surface $\Sigma'''$ arising in our construction satisfies $\mathcal{C}^{2, \alpha}$-estimates depending only on $\mathcal{H}^{n-1}(\partial_1 \Omega)$, $|p|_{\mathcal{C}^1(\bar\Omega)}$, and the local geometry of $(M^n, g)$. It is stable in the sense of Appendix \ref{sec: remark on SS81}, it minimizes area with respect to variations in $\bar{\Omega}'''$, and it is $C$-almost minimizing in $\Omega'\cap\Omega''$. 

\begin{proof}
Let $\{\Sigma'''_\varepsilon\}$ be the surfaces constructed in Lemma \ref{lem: kicking out the kink}. Using the uniform volume- and $\mathcal{C}^{2, \alpha}$-estimates (alternatively using Theorem \ref{thm: varifold compactness for stable submanifolds}), we can pass these surfaces to an immersed subsequential limit $\Sigma^{n-1}$. By the $C$-almost minimizing property of $\Sigma_\varepsilon^{n-1}$ in $\Omega' \cap \Omega''$ which descends to $\Sigma^{n-1}$, it follows that $\Sigma^{n-1}$ could only fail to be embedded on the boundary of $\Omega' \cap \Omega''$. As above, the argument in Corollary \ref{cor: increasing limits are embedded} rules out sheeting near this boundary. Hence $\Sigma^{n-1}$ is an embedded generalized apparent horizon. Finally, the strong maximum principle applied as in the proof of Theorem \ref{thm: existence of generalized apparent horizon in Omega} shows that if a component of $\Sigma^{n-1}$ touches a component of $\Sigma'$ or $\Sigma''$ then these components must coincide.
\end{proof}
\end{proposition}

\noindent We emphasize that $\Sigma'''$ in this proposition is stable and outer minimizing. Hence the conclusion is not trivial when $\Sigma' \cap \Sigma'' = \emptyset$ or even when $\Omega' = \Omega''$. This will be relevant in the proof of Theorem  \ref{thm: existence of outermost generalized apparent horizon}. 

\begin{remark}
The particular $\varepsilon$-regularization we chose helped us with the construction of appropriate super solutions for the problems $(\mc - |\tr(p)|_\varepsilon - t) u_t^\varepsilon = 0$ from the conditions $\mc_{\Sigma'} = |\tr_{\Sigma'}(p)|$ and $\mc_{\Sigma''} = |\tr_{\Sigma''}(p)|$ in Lemma \ref{lem: kicking out the kink}. For marginally outer trapped surfaces, where the natural conditions for the inner boundaries are $\mc_{\Sigma'} + \tr_{\Sigma'}(p) \leq 0$ and $\mc_{\Sigma''} + \tr_{\Sigma''}(p) \leq 0$, one can proceed similarly by first finding auxiliary surfaces $\Sigma^{n-1}_\varepsilon \subset \Omega' \cap \Omega''$ for which $\mc_{\Sigma_\varepsilon} + \tr_{\Sigma_\varepsilon} p_\varepsilon = 0$. Here, $p_\varepsilon := p - \varepsilon \phi g$ for some fixed smooth function $\phi: \Omega \to \mathbb{R}$ with $\phi \equiv 1$ near $\Sigma' \cup \Sigma''$ and $\phi \equiv 0$ near the outer boundary $\partial_1 \Omega$. The apparent horizon enclosing $\Sigma' \cup \Sigma''$ is then found by letting $\varepsilon \searrow 0$ as above (note that the $C$-almost minimizing property is independent of $\varepsilon>0$ here as well, so we have geometric estimates that allow us to pass to a limit). The trick of modifying the second fundamental form tensor in the direction of $g$ to get strict barriers was used in \cite[Theorem 5.1]{AM07}. 
\end{remark}


\section{The outermost Generalized Trapped Surface} \label{sec: existence of outermost generalized apparent horizon}

In this section we give the proof of Conjecture \ref{conj: Ilmanen}. The outermost generalized trapped surface is constructed as the boundary of the union of all generalized trapped domains. This is analogous to the construction of the apparent horizon as the boundary of the trapped region in \cite{HawEll}, \cite{KrieHay}, \cite{HuiIlm}, \cite{AM07}. The idea of replacing the total union by one increasing union is contained in \cite{KrieHay} and has been used in \cite{AM07} to prove existence of a smooth outermost apparent horizon in $3$-dimensional data sets. A major technical challenge in \cite{AM07} was to show that the boundary of this increasing union is smooth and embedded. We survey the steps in their argument for comparison: the authors first reduce to the case where the boundaries of the increasing sets in this union are stable marginally outer trapped surfaces and hence have bounded curvature \cite[Theorem 1.2] {AM05} (see also \S4 in \cite{HuiIlm}), and by a further reduction to the case where these surfaces have a lower bound on their ``outer injectivity radius," cf. \cite [\S 6] {AM07}. The latter step depends on the a priori curvature bound coming from stability and a very delicate surgery procedure. The lower bound on the outer injectivity radius and the curvature bound give an area estimate \cite[Theorem 6.1]{AM07} for the boundary surfaces in this union. It follows that the limit of these surfaces exists as a smooth immersed marginally trapped surface which is the boundary of an open set, and hence cannot touch itself on the \emph{inside}. The lower bound on the \emph{outer} injectivity radius of these surfaces implies that the limit is embedded \cite[\S 7]{AM07}. \\ 

\noindent Our proof of Conjecture \ref{conj: Ilmanen}, which also applies to marginally outer trapped surfaces and works (at least) in dimensions $2 \leq n \leq 7$, is based on the lower order properties (outward minimizing and outward almost minimizing) of these surfaces:

\begin{theorem} \label{thm: existence of outermost generalized apparent horizon} Let $(M^n, g, p)$ be a complete initial data set and assume that $2 \leq n \leq 7$. Let $\Omega \subset M^n$ be a bounded open subset with embedded $\mathcal{C}^2$-boundary $\partial \Omega = \partial_1 \Omega \dot{\cup} \partial_2 \Omega$ such that $\mc_{\partial_1 \Omega} > |\tr_{\partial_1 \Omega}(p)|$ (mean curvature with respect to the normal pointing out of $\Omega$) and $\mc_{\partial_2 \Omega} \leq |\tr_{\partial_2 \Omega}(p)|$ (with respect to the normal pointing into $\Omega$). Then there exists $\Omega' \subset \Omega$ open such that $\partial \Omega' = \partial_1 \Omega \dot \cup \Sigma'$ where $\Sigma'$ is an embedded $\mathcal{C}^{2, \alpha}$-hypersurface with $\mc_{\Sigma'} = |\tr_{\Sigma'}(p)|$ (with respect to the normal pointing into $\Omega'$), and such that $\Sigma'$ is the outermost generalized trapped surface in the following sense: if $\Omega'' \subset \Omega$ has boundary $\partial_1 \Omega \dot \cup \Sigma''$ where $\Sigma''$ is an embedded $\mathcal{C}^{2}$-hypersurface with $\mc_{\Sigma''} \leq |\tr_{\Sigma''}(p)|$, then $\Omega' \subset \Omega''$ (so that $\Sigma'$ encloses $\Sigma''$). Moreover, this $\Sigma'$ minimizes area with respect to variations in $\bar {\Omega}'$.  

\begin{proof}

Consider the closed set $F = \cap \bar{\Omega}''$ where the intersection is taken over all open subsets $\Omega'' \subset \Omega$ for which $\partial \Omega'' = \partial_1 \Omega \dot \cup \Sigma''$ such that $\Sigma''$ is an embedded, $\mathcal{C}^{2}$-hypersurface with $\mc_{\Sigma''} \leq |\tr_{\Sigma''}(p)|$. As in \cite{KrieHay}, $F$ is already the intersection of a countable family $\bar \Omega''_i$ of such sets. (This is because $\Omega$ is second countable.) We can use Proposition  \ref{prop: to any two generalized apparent horizons that intersect there is another one enclosing them} to arrange for these sets to be decreasing $ \Omega''_1 \supset \Omega''_2 \supset \ldots $ and such that $\Sigma''_i$ ($= \partial \Omega''_i \setminus \partial_1 \Omega$) is a stable (in the sense of inequality (\ref{inequ: stability of graph(u)}) and Appendix \ref{sec: remark on SS81}) generalized apparent horizon (cf. \cite{KrieHay}, \cite{HuiIlm}, \cite{AM07}) which minimizes area with respect to variations in $\bar {\Omega}''_i$. Then Corollary \ref{cor: increasing limits are embedded} implies that $\Omega' := \interior F$ has the required properties. 
\end{proof}
\end{theorem}

\appendix


\section{A remark on \cite{SchoenSimon}} \label{sec: remark on SS81}

In this appendix we explain how the regularity theory of \cite{SchoenSimon} for stable critical points of elliptic variational problems can be applied to the context of generalized apparent horizons, even though these latter surfaces are not associated with a particular functional. The proofs in \cite{SchoenSimon} generalize to our situation with only a few very minor modifications, which we discuss here. The results stated in this appendix provide a general $\mathcal{C}^{1, \alpha}$-regularity and compactness theory for limits of smooth embedded hypersurfaces of bounded area and mean curvature that also satisfy a stability-type inequality  such as (\ref{inequ: stability of graph(u)}), with the usual estimate of the singular set of such a limit. In particular, these results are available in all dimensions and in other situations where a priori curvature estimates may not be readily available. \\

\noindent We follow the notation and conventions of \cite{SchoenSimon} closely in this appendix to facilitate reference. \\

Given a Riemannian manifold $(M^{n}, g)$, a point $p \in M^n$, and $0 < \rho_0 < \inj_p(M^n, g)$, one can use geodesic normal coordinates centered at $p$ to identify the geodesic ball $B^{n}(p, \rho_0) \subset M^{n}$ with the Euclidean ball $\{|X| < \rho_0\} \subset \Tan_p(M^n)$. Given a hypersurface $G^{n-1} \subset M^n \cap B^{n}(p, \rho_0)$ one can use this identification to compute geometric quantities of $G^{n-1}$ either with respect to $g_{ij}$ or with respect to the Euclidean metric $\delta_{ij}$ on $\{|X| < \rho_0\}$. One has that
\begin{eqnarray*}
g_{i j}\big(0) &=& \delta_{i j}, \\ 
\partial_k g_{i j}\big(0) &=& 0,  \text{ and} \\
\sup_{\{|X| < \rho_0\}} \left|\partial^2_{k l} g_{ij} \right| &\leq& \mu_1^2
\end{eqnarray*}
for some constant $\mu_1 \geq 0$ depending only on the geometry of $(B^n(p, \rho_0), g)$ (specifically on the $\mathcal{C}^0$ size of the curvature tensor). Denoting all quantities computed with respect to the Euclidean metric $\delta_{ij}$ with a hat, one has that 
\begin{eqnarray} \label{eqn: compare g metric to delta metric}
\left| |\h|_g^2  - |\hat{\h}|^2 \right| &\leq& c_1 \left( \mu_1 |X| |\hat{\h}|^2 + \mu_1^2\right) \\
\left| \mc - \hat \mc \right| &\leq& c_1 \left(\mu_1 |X| |\hat{\h}| + \mu_1 \right) \nonumber
\end{eqnarray}
\noindent provided that $\mu_1 |X| < 1$. Here, $X$ is the position vector in $\{ |X| < \rho_0 \}$, $c_1$ is a dimensional constant, and the norms of hatted quantities are taken with respect to the Euclidean metric $\delta_{ij}$. Cf. \cite[\S 1, \S 6]{SchoenSimon}. \\

Suppose now that the embedded $\mathcal{C}^2$-hypersurface $G^{n-1} \subset M^n$ satisfies a stability-type inequality of the form 
\begin{equation} \label{appendix: stability of graph(u)}
\left(1 - \frac{\eta}{2} \right) \int_{G} |\h|_g^2 \phi^2 \leq \int_G |\bar{\D}_g \phi|^2  + \kappa^2 \int_G \phi^2 \ \ \forall \phi \in \mathcal{C}^1_c(M^n)
\end{equation}
where $\eta \in (0, 1)$ and that in addition its mean curvature is bounded $|\mc| \leq \kappa$. Using the estimates in (\ref{eqn: compare g metric to delta metric}) we see that on $B^n(p, \rho_0)$ this inequality carries over to the Euclidean geometry in the form  
\begin{eqnarray} \label{appendix: stability of graph(u) in Euclidean}
\left(1 - \eta\right) \int_G |\hat{\h}|^2 \phi^2 \leq \\ \nonumber \int_G |\hat{\D} \phi|^2  + \hat{\kappa}^2 \int_G \phi^2 \ \ \text{ for all } \phi \in \mathcal{C}^1_c(\{|X| < \rho_0 \})
\end{eqnarray}
\noindent provided that $\mu_1 \rho_0$ is sufficiently small (depending only on the dimension), and where the new constant $\hat{\kappa}$ depends on $\kappa, \mu_1, \eta$ and the dimension $n$. All integrals here are computed with respect to the Euclidean metric induced on $G$, and $\hat{\D} \phi$ denotes the tangential gradient \emph{along} $G$. From (\ref{eqn: compare g metric to delta metric}) we see that the mean curvature $\hat{\mc}$ of $G$ on the geodesic ball $\{|X| < \rho_0 \}$ can be estimated by
\begin{equation} \label{eq: estimate on mean curvature}
|\hat {\mc}| \leq c_1 \left(\mu_1 |X||\hat \h| + \mu_1\right) + \kappa.
\end{equation}
\noindent Absorbing $\kappa$ into the constant $\mu_1$, we see that the structural assumptions $(1.16)$ and $(1.17)$ of \cite{SchoenSimon} are satisfied with a marginally worse constant $1 - \eta$ multiplying the left-hand side of the stability-type inequality (\ref{appendix: stability of graph(u) in Euclidean}). We note that the fundamental integral curvature estimate in \cite[\S2]{SchoenSimon} still follows from these inequalities provided that $0 < \eta < \frac{1}{n-1}$ and that the proof of the basic regularity estimate \cite[Theorem 1]{SchoenSimon} carries over \emph{verbatim} to our setting if we assume a priori that $G^{n-1} \subset B^n(p, \rho_0)$ is an embedded, relatively closed $\mathcal{C}^{2}$-hypersurface. In the statement of the result below, $C(X, \rho)$ denotes the cylinder $\{y \in \mathbb{R}^{n-1} : |y - x | < \rho \} \times \mathbb{R}$ for  $X = (x, x_{n+1}) \in \{ |X| < \rho_0 \}$: 

\begin{theorem} [\cite{SchoenSimon}] \label{thm: basic regularity in SS81} Let $G^{n-1} \subset \{|X| < \rho_0 \}$ be a relatively closed, embedded $\mathcal{C}^{2}$-hypersurface such that for some constants $\mu, \mu_1 > 0$ one has 
\begin{eqnarray*}
&&\mathcal{H}^{n-1}(G^{n-1}) \leq \mu \omega_{n-1} \rho_0^{n-1},  \\
&&|\hat {\mc}| \leq c_1 \left(\mu_1|X||\hat \h| + \mu_1 \right),
\end{eqnarray*} and
\begin{eqnarray*}
\left(1 - \eta \right) \int_G |\hat{\h}|^2 \phi^2 \leq \int_G |\hat{\D} \phi|^2  + \mu_1^2 \int_G \phi^2  \text{ for all } \phi \in \mathcal{C}^1_c(\{|X| < \rho_0 \})
\end{eqnarray*}
for some $\eta \in (0, \frac{1}{n-1})$. Then there exists a number $\delta_0 \in (0, 1)$ depending only on $n, \mu, \eta$, and  $\mu_1 \rho_0$ such that if $X \in G^{n-1} \cap \{|X| < \rho_0/4 \}$, $\rho \in (0, \rho_0/4)$, and $G'$ is the connected component of $G^{n-1} \cap C(X, \rho)$ containing $X$, and if for some $\delta \in (0, \delta_0)$
\begin{eqnarray*}
\sup_{Y \in G'} |y_{n+1} - x_{n+1}| &\leq& \delta \rho \\
\mu_1 \rho &\leq& \delta
\end{eqnarray*}
where $X = (x, x_{n+1})$ and $Y = (y, y_{n+1})$, then $G' \cap C(X, \rho/2)$ consists of a disjoint union of graphs of functions $u_1 < u_2 < \ldots < u_k$ defined on $B^{n-1}(x, \rho/2) \subset \mathbb{R}^{n-1}$ such that
\begin{eqnarray} \label{est: basic regularity in SS81}
\frac{2}{\rho} \sup_{y \in B^{n-1}(x, \rho/2)}|u_i (y)| + \sup_{y \in B^{n-1}(x, \rho/2)} |\hat \D u_i (y)| + \\ \nonumber \left(\frac{\rho}{2}\right)^{\alpha} \sup_{y, y' \in B^{n-1}(x, \rho/2), y \neq y'} \frac{|\hat \D u_i (y) - \hat \D u_i (y')|}{|y - y'|^\alpha} =:|u_i|_{1, \alpha, B^{n-1}(x, \rho/2)} &\leq& \delta^{\frac{1}{3}}
\end{eqnarray}
for all $i = 1, 2, \ldots k$, where $\alpha \in (0, 1)$, $k$, and $c_2$ depend only on $n, \mu, \eta$, and  $\mu_1 \rho_0$. 
\end{theorem}

\begin{remark} \label{rem: basic regularity in SS81} The appearance of $\delta^{\frac{1}{3}}$ in estimate (\ref{est: basic regularity in SS81}) rather than plain $\delta$ as in \cite{SchoenSimon} is due to the fact that we don't want to appeal to Schauder theory at this point (to avoid mention of a defining equation for the functions $u_i$). Careful screening of the proof in \cite{SchoenSimon} (note in particular (1.21), (3.32), (4.33)-(4.37) in that paper) shows that this power is sufficient. Theorem \ref{thm: basic regularity in SS81} says that a regular closed embedded submanifold with controlled mean curvature and area that is stable in the sense of (\ref{appendix: stability of graph(u) in Euclidean}) decomposes into a union of graphs with $\mathcal{C}^{1, \alpha}$-estimates whenever its support is sufficiently close to a hyperplane. In the case where the submanifolds are stationary with respect to an elliptic functional (i.e. satisfy an appropriate equation), Theorem \ref{thm: basic regularity in SS81} is combined in \cite{SchoenSimon} with the compactness theorem for rectifiable varifolds, a version of Federer's dimension reduction argument, and the fact that there exist no stable minimal hypercones in $\mathbb{R}^n$ other than planes when $3 \leq n \leq 7$ to obtain curvature estimates in these dimensions. (Note that there are stable minimal tangent cones in dimension $n=2$ which are singular at the origin. Such tangent cones for the limiting surfaces are ruled out in \cite[pages 786 and 787]{SchoenSimon}. We point out that in this paper, one-dimensional generalized apparent horizons arise as cross-sections of stable almost minimizing cylinders and that we can carry out the regularity argument on that level. The almost minimizing property and the fact that all our surfaces are boundaries would also rule out such singular cones.) If we do bring in the defining equation $\mc_G = |\tr_G (p)|_\varepsilon$ of the surfaces in this paper, uniform $\mathcal{C}^{2, \alpha}$-estimates in terms of $\mu, \mu_1 \rho_0, n$ and $|p|_{\mathcal{C}^1}$ follow from the same argument. 
\end{remark}

The following regularity property was formulated in \cite[page 780]{SchoenSimon}: 

\begin{definition} [\cite{SchoenSimon}] \label{def: property P} Let $\delta_0, \mu, \mu_1$ be positive constants. We say that a countably $(n-1)$-rectifiable varifold $V \in \mathbf{V}_{n-1}(\mathbb{R}^n)$ has the property $\mathbf{P}_{\delta_0 \mu \mu_1}$ with respect to an open subset $U$ of $\{|X| < \rho_0 \}$ if
\begin{eqnarray} \label{eq: upper and lower bound for density}
\mu^{-1} \sigma^{n-1} \omega_{n-1} \leq ||V||(B^n(X, \sigma)) \leq \mu \sigma^{n-1} \omega_{n-1} \\ \nonumber \text{ for all } X \in \spt ||V|| \text{ and } \sigma > 0 \text{ such that } B^n(X, \sigma) \subset U 
\end{eqnarray}
and provided that whenever the hypotheses
$Y \in \spt ||V||, B^n(Y, \rho) \subset U, \mu_1 \rho < \delta$ and $\spt ||V|| \cap B^n(Y, \rho) \subset \{X \in \mathbb{R}^n : \dist(X, \Pi) < \delta \rho \}$ hold for some hyperplane $\Pi \subset \mathbb{R}^n$ containing $Y$ and some $\delta \in (0, \delta_0)$, then there exists an isometry $\mathcal{O}$ of $\mathbb{R}^n$ with $\mathcal{O}(Y) = 0$, $\mathcal{O}(\Pi) = \mathbb{R}^{n-1} \times \{0\}$, and 
\begin{equation*}
\mathcal{O}(\spt ||V|| \cap B^n(Y, \rho)) \cap C(0, \rho/2) = \bigcup_{i=1}^l \graph(u_i)
\end{equation*}
for some integer $l$, $1 \leq l \leq \mu$, where $u_i \in \mathcal{C}^{1, \alpha}(B^{n-1}(0, \rho/2))$ are such that $u_1 \leq u_2 \leq \ldots \leq u_l$ and
\begin{equation*}
|u_i|_{1, \alpha, B^{n-1}(0, \rho/2)} \leq {\delta}^{\frac{1}{3}}. 
\end{equation*}
Here we are using scale-invariant Schauder norms on the left. 
\end{definition}

\begin{remark} The upper bound in condition (\ref{eq: upper and lower bound for density}) is implied by a total mass bound $||V||(U)$ and an estimate on the first variation $\delta V$ of $V$ of the form \begin{eqnarray} \label{eqn: bounded mean curvature} \big| \delta V ( \varphi ) | \leq \mu_1 \int \left( |\varphi| + |X| |\D \varphi| \right) d ||V|| \text{ for all } \varphi \in \mathcal{C}^1_c (U, \mathbb{R}^n) \end{eqnarray} through the monotonicity formula, see \cite[pages 778, 779]{SchoenSimon}. Here, $\D \varphi$ denotes the ambient covariant derivative of the vector field $\varphi$. Such an estimate is implied by (\ref{eq: estimate on mean curvature}). The lower bound in  (\ref{eq: upper and lower bound for density}) follows from the same principle if we assume that $V$ satisfies a positive lower bound on its density (for example if it is integer multiplicity). Note that we \emph{do not} require that the defining graphs $u_1 \leq \ldots \leq u_l$ be disjoint. It is also evident that the class $\mathbf{P}_{\delta_0 \mu \mu_1}$ is preserved under varifold limits (by Arzela-Ascoli). Note also that any relatively closed, embedded $\mathcal{C}^2$-hypersurface $G^{n-1} \subset \{|X| < \rho_0 \}$ satisfying the hypotheses of Theorem \ref{thm: basic regularity in SS81} belongs to the class $\mathbf{P}_{\delta_0 \mu \mu_1}$ provided one chooses $\delta_0$ sufficiently small, cf. \cite[Remark 11]{SchoenSimon}. 
\end{remark}

The following definitions of the singular and regular sets differ marginally from the definition in \cite[page 777]{SchoenSimon}. 

\begin{definition} Let $U \subset \{|X| < \rho_0 \}$ be an open set and let $V \in \mathbf{V}_{n-1}(\mathbb{R}^n)$ have the property $\mathbf{P}_{\delta_0 \mu \mu_1}$ with respect to this set. The regular set $\reg (V)$ of the varifold $V$ in $U$ is defined as the set of all points $Y \in \spt ||V|| \cap U$ such that for some small radius $\rho>0$ one can write $\spt ||V|| \cap B(Y, \rho)$ as the union of weakly ordered $\mathcal{C}^{1, \alpha}$-graphs $u_1 \leq u_2 \leq \ldots \leq u_l$ defined on a common hyperplane. The singular set $\sing(V)$ is defined as the complement of $\reg(V)$ in $\spt(V) \cap U$. 
\end{definition}

\noindent Note that it follows from the property $\mathbf{P}_{\delta_0 \mu \mu_1}$ that $Y \in \reg(T)$ if and only if there exists a varifold tangent $T \in \Tan \Var (V, Y)$ so that $||T|| = m |\Pi|$ for some hyperplane $\Pi$ of $\mathbb{R}^n$ and some positive integer $0 < m \leq \mu$. It  also follows trivially that varifold tangents at regular points are unique. From the definition one sees that $\reg(V)$ is relatively open. There is some minor subtlety in allowing the graphs of the locally defining functions $u_1 \leq u_2 \leq \ldots \leq u_k$ to touch. This is related to the fact that, unlike minimal surfaces, (generalized) apparent horizons don't satisfy a two-sided maximum principle. If the first variation of $V$ also satisfies (\ref{eqn: bounded mean curvature}), then its tangent varifolds are stationary cones \cite[page 780]{SchoenSimon}. The Hopf maximum principle and the constancy theorem then show that the tangent cones of such varifolds are smooth hypersurfaces with constant (integer) multiplicity near their regular points. \\

\noindent It is easy to see that this notion of singular set is consistent with the basic assumptions A.1, A.2, and A.3 of the abstract dimension reduction procedure given in Appendix A of \cite{LeonGMT}. The compactness theory of \cite{SchoenSimon} for stable minimal hypersurfaces with area bounds takes the following form in the present context, with virtually the same proof: 

\begin{theorem} [\cite{SchoenSimon}] \label{thm: varifold compactness for stable submanifolds} Let $G^{n-1}_i \subset \{|X| < \rho_0 \}$ be a sequence of embedded relatively closed $\mathcal{C}^2$-hypersurfaces such that there exist constants $\mu_1, \mu > 0$ and $\eta \in (0, \frac{1}{n-1})$ so that
\begin{eqnarray*}
&&\mathcal{H}^{n-1}(G^{n-1}_i) \leq \mu \rho_0^{n-1} \omega_{n-1} \\
&&\left(1 - \eta \right) \int_{G_i} |\hat{\h}_i|^2 \phi^2 \leq \int_{G_i} |\hat{\D} \phi|^2  + \mu_1^2 \int_{G_i} \phi^2 \ \forall \phi \in \mathcal{C}^1_c(\{|X| < \rho_0 \})\\ 
&&|\hat {\mc}_i| \leq c_1 \left(\mu_1 |X||\hat{\h}_i| + \mu_1 \right). 
\end{eqnarray*}
Let $G^{n-1} \in \mathbf{V}_{n-1}(\mathbb{R}^n)$ be a subsequential varifold limit of $\{G_i^{n-1}\}$. Then $G^{n-1}$ is integer rectifiable and the Hausdorff dimension of its singular set $\sing(G^{n-1})$ is $\leq n-8$. In particular, if $2 \leq n \leq 7$, the limit $G^{n-1}$ is an immersed $\mathcal{C}^{1, \alpha}$-hypersurface. 
\end{theorem}

In the following corollary we note that increasing limits of stable outer minimizing hypersurfaces with uniformly bounded mean curvature remain regular and embedded.   

\begin{corollary} \label{cor: increasing limits are embedded} Assume that $2 \leq n \leq 7$, let $(M^n, g)$ be a complete Riemannian manifold, and let $\Omega \subset M^n$ be a bounded open set with smooth boundary $\partial \Omega = \partial_1 \Omega \dot \cup \partial_2 \Omega$. Consider a decreasing sequence $\Omega \supset \Omega_1 \supset \Omega_2 \supset \ldots$ of open subsets of $\Omega$ with $\partial \Omega_i = \partial_1 \Omega  \dot \cup \Sigma^{n-1}_i$ such that the $\Sigma^{n-1}_i$ are embedded $\mathcal{C}^{2}$-hypersurfaces which are outer minimizing, i.e., minimize area with respect to variations in $\bar \Omega_i$. Assume that there are constants $\mu_1, \mu$ so that the assumptions of Theorem \ref{thm: varifold compactness for stable submanifolds} hold uniformly for $\Sigma_i^{n-1}$ and that the surfaces $\Sigma^{n-1}_i$ stay away in Hausdorff distance from the outer boundary $\partial_1 \Omega$. Then there exists an open set $\Omega' \subset \Omega$ with $\partial \Omega' =  \partial_1 \Omega \dot \cup \Sigma^{n-1}$ such that $\Sigma^{n-1}$ is an embedded $\mathcal{C}^{1, \alpha}$-hypersurface, $\Sigma_i^{n-1} \to \Sigma^{n-1}$ in $\mathcal{C}^{1, \beta}$ for any $0 < \beta < \alpha$, and so that $\Sigma^{n-1}$ minimizes area with respect to variations in $\bar \Omega'$.
\end{corollary}

\noindent The corollary holds true if `outer minimizing' is replaced by `outer $C$-almost minimizing' with essentially the same proof. The uniform area bound (expressed in the constant $\mu$) required for the use of Theorem \ref{thm: varifold compactness for stable submanifolds} is given by $\mathcal{H}^{n-1}(\partial_1 \Omega)$, respectively  $\mathcal{H}^{n-1} (\partial_1 \Omega) + C \mathcal{L}^n(\Omega)$. The H\"older exponent $\alpha \in (0, 1)$ is as in Theorem \ref{thm: basic regularity in SS81} and depends only on $n$, $\eta$, $\mu$, and $\mu_1 \rho_0$. 

\begin{proof} Note first that by Allard's integral compactness theorem, the sequence $\Sigma^{n-1}_i$ converges to a countably $(n-1)$-rectifiable integer multiplicity varifold $\Sigma^{n-1}$ with bounded mass and first variation. From Theorem \ref{thm: varifold compactness for stable submanifolds} we know that $\sing(\Sigma^{n-1}) = \emptyset$. Let $T = m|\Pi|$ be the (unique) varifold tangent at $X \in \spt \Sigma^{n-1}$ where $\Pi$ is a hyperplane and $m$ is a positive integer. The bounds on the mean curvature imply that $T \cap B^n(0, 1)$ is approached in Hausdorff distance by appropriate rescalings of the embedded hypersurfaces $\Sigma_i^{n-1}$ (cf. \cite[page 780]{SchoenSimon}), which by Theorem \ref{thm: basic regularity in SS81} decompose after an appropriate rotation into graphs $u_1^i < u_2^i < \ldots < u_m^i$ over $B^{n-1}(0, 1/2) \subset \mathbb{R}^{n-1}$ with $|u|_{1, \alpha, B^{n-1}(0, 1/2)} \to 0$. Since the $\Sigma_i^{n-1}$ are boundaries, the set of points $\{(x, x_n) : x \in B^{n-1}(0, 1/2) \text{ and } u^i_{j}(x) < x_n < u^i_{j+1}(x) \}$ either belongs to $\Omega_i$ or its complement for every $j = 1, \ldots, m-1$. Since $\Sigma^{n-1}_i$ minimizes area in $\bar \Omega_i$ we immediately obtain that $m \leq 2$. Finally note that since the sets $\Omega_i$ are decreasing, $m = 1$ (cf. \cite[page 971]{AM07}). Since $X \in \spt \Sigma^{n-1}$ was arbitrary it follows that $\Sigma^{n-1}$ is embedded. That $\Sigma^{n-1} = \partial \Omega'$ and that $\Sigma^{n-1}$ minimizes area with respect to variations in $\bar \Omega '$ now follow easily. 
\end{proof}

\begin{remark} \label{rem: curvature estimates follow too} Theorem \ref{thm: varifold compactness for stable submanifolds} describes the regularity of varifold limits of smooth embedded hypersurfaces which have bounded area and mean curvature, and which satisfy a uniform stability-type inequality. An immediate consequence are curvature estimates for stable embedded generalized apparent horizons $\Sigma^{n-1}$ in dimensions $2 \leq n \leq 7$ in terms of the injectivity radius of $(M^n, g)$, $|\Rm|_{\infty}$ (which enters through $\mu_1$), $|p|_{\mathcal{C}^1(\bar \Omega)}$, and a bound on the their hypersurface measure (coming from the outer minimizing property of the surfaces in this paper). Such estimates have been obtained for immersed stable marginally outer trapped surfaces in dimensions $2 \leq n \leq 6$ in \cite{AM05}, by generalizing the iteration method in \cite{PMTII}, in particular by deriving an appropriate analogue of the Simons identity for $\chi$ (which is the second fundamental form plus the restriction of $p$ to the surface). The authors derive an inequality like (\ref{appendix: stability of graph(u)}) with $|h|$ replaced by $|\chi|$ using the first eigenfunction of the stability operator associated to a marginally trapped surface in an initial data set in \cite[Lemma 5.6]{AM05}. Their lemma shows that the theory in this appendix also applies to stable marginally trapped surfaces. The `lower order approach' of \cite{SchoenSimon} that we are taking here is quite flexible and applies nicely to the class of generalized apparent horizons for which one expects $\mathcal{C}^{2, \alpha}$-regularity at best. These results are available in all dimensions if we accept singular sets of Hausdorff co-dimension $7$. Note that the important estimate on the outer injectivity radius in \cite[\S6]{AM07} can be recovered for stable outer $C$-almost minimizing surfaces $\Sigma^{n-1}$ in all dimensions $2 \leq n \leq 7$ by combining Theorem \ref{thm: varifold compactness for stable submanifolds} with the argument of Corollary \ref{cor: increasing limits are embedded}. The point is that the increasing property of the surfaces $\Sigma^{n-1}_i$ in Corollary \ref{cor: increasing limits are embedded} only enters when we concluded that $m=1$. That there can be at most two sheets merging from the inside follows from the one-sided (almost) minimizing property and the fact that the $\Sigma^{n-1}_i$'s are boundaries. 
\end{remark}

\end{document}